\documentclass{article}
\usepackage{latexsym,amssymb}
\input amssym.def \input amssym

\newtheorem{te}{Theorem}[section]

\newtheorem{de}[te]{Definition}
\newtheorem{lm}[te]{Lemma}
\newtheorem{fa}[te]{Fact}
\newtheorem{pp}[te]{Proposition}

\newtheorem{ex}[te]{Example}
\newtheorem{qu}[te]{Question}
\def\dokaz{\noindent{\bf Proof. }}
\def\kraj{\hfill $\Box$ \par \vspace*{2mm} }
\newcommand{\zve}[1]{{{}^*\hspace{-0.5mm}#1}}
\newcommand{\zvez}[1]{{{}^*\hspace{-1mm}#1}}

\def\widemid{\hspace{1mm}\widetilde{\mid}\hspace{1mm}}
\def\zvezmid{\hspace{1mm}\zvez\mid\hspace{1mm}}
\def\zvepar{\hspace{1mm}\zvez\parallel\hspace{1mm}}
\def\nwidemid{\hspace{1mm}\widetilde{\nmid}\hspace{1mm}}
\def\nzvemid{\hspace{1mm}\zvez\nmid\hspace{1mm}}
\def\lev{{\rm lev}}
\def\zs{\forall}
\def\po{\exists}
\def\str{\rightarrow}
\def\Str{\Rightarrow}

\def\dl{\Leftrightarrow}

\def\gstr{\hspace{-0.1cm}\uparrow}
\def\ps{\subseteq}

\def\dom{{\rm dom}}
\def\ran{{\rm ran}}

\def\rest{\upharpoonright}

\def\cF{{\cal F}}
\def\cG{{\cal G}}
\def\cH{{\cal H}}
\def\cU{{\cal U}}

\def\cW{{\cal W}}

\def\cP{{\cal P}}
\def\cQ{{\cal Q}}

\begin{document}
\begin{center}
           {\huge \bf Divisibility in $\beta N$ and $\zve N$}
\end{center}
\begin{center}
{\small \bf Boris  \v Sobot}\\[2mm]
{\small  Department of Mathematics and Informatics, Faculty of Sciences,\\
University of Novi Sad,\\
Trg Dositeja Obradovi\'ca 4, 21000 Novi Sad, Serbia\\
e-mail: sobot@dmi.uns.ac.rs\\
ORCID: 0000-0002-4848-0678}
\end{center}
\begin{abstract} \noindent
The paper first covers several properties of the extension of the divisibility relation to a set $\zve N$ of nonstandard integers, including an analogue of the basic theorem of arithmetic. After that, a connection is established with the divisibility in the Stone-\v Cech compactification $\beta N$, proving that the divisibility of ultrafilters introduced by the author is equivalent to divisibility of some elements belonging to their respective monads in an enlargement. Some earlier results on ultrafilters on lower levels on the divisibility hierarchy are illuminated by nonstandard methods. Using limits by ultrafilters we obtain results on ultrafilters above these finite levels, showing that for them a distribution by levels is not possible.
\vspace{1mm}\\

\end{abstract}

In \cite{So1} four different ways to extend the divisibility relation on $N$ to the Stone-\v Cech compactification  $\beta N$ were introduced. One of them, the relation $\widemid$, seems to be the most fruitful for investigation, and some of its properties were extracted in \cite{So3}. In this paper we prove more about it, using mostly the connection between $\beta N$ and nonstandard extensions $\zve N$ of $N$.

The first section covers most of the facts needed for the rest of the paper about nonstandard extensions, $\beta N$ and the connection between them. In the second section we single out some facts on the divisibility in $\zve N$ that can be useful for obtaining information on divisibility in $\beta N$. The new results are mostly contained in sections 3 and 4, the former containing facts on the first $\omega$-many levels of the $\widemid$-hierarchy and the latter on ultrafilters above these finite levels. Some open problems are addressed in the last section.

\section{Introduction}

{\bf Nonstandard methods.} In the course of the last 60 years many approaches to nonstandard extensions have been developed. One general idea is to work with a particular construction of nonstandard universe, most frequently an ultrapower or an ultralimit. The other is to give an axiomatic development of the notion of a nonstandard extension. The paper \cite{N1} contains an overview of various historically relevant axiomatic systems. In this paper we mostly follow the Robinson-Zakon superstructure approach as exposed in Henson's chapter \cite{H}, a text highly recommended for mathematicians new to the subject.\\

Let $X$ be a set; we assume that elements of $X$ are atoms: none of them contains as an element any of the others. Let $V_0(X)=X$, $V_{n+1}(X)=V_n(X)\cup P(V_n(X))$ for $n\in\omega$ and $V(X)=\bigcup_{n<\omega}V_n(X)$. The rank of $x\in V(X)$ is the smallest $n\in\omega$ such that $x\in V_n(X)$. $V(X)$ is called a {\it superstructure}. We call atoms (elements of $X$) and sets of $V(X)$ by a common name {\it objects} in $V(X)$.

Let $V(X)$ be a superstructure. Its {\it nonstandard extension} is a pair $(V(Y),*)$, where $V(Y)$ is a superstructure with the set of atoms $Y\supset X$ and $*:V(X)\str V(Y)$ is a rank-preserving function such that $\zve X=Y$ and satisfying the following principle.\\

{\it The Transfer Principle.} For every bounded formula $\varphi$ and every $a_1,a_2,\dots$, $a_n\in V(X)$, $\varphi(a_1,a_2,\dots,a_n)$ holds in $V(X)$ if and only if $\varphi(\zve a_1,\zve a_2,\dots,\zve a_n)$ holds in $V(Y)$.\\

(A first-order formula is bounded if all its quantifiers are bounded, i.e.\ of the form $(\zs x\in y)$ or $(\po x\in y)$. The free variables that appear in $\varphi(a_1,a_2,\dots,a_n)$ are exactly objects $a_1,a_2,\dots,a_n$ from $V(X)$ and in $\varphi(\zve a_1,\zve a_2,\dots,\zve a_n)$ they are replaced with their star-counterparts. The atomic subformulas in $\varphi$ are of the form $A(x_1,\dots,x_k)$ for some $k$-ary relation $A\in V(X)$.) 

We may abuse the notation and call $V(\zve X)$ the extension of $V(X)$, or even call $\zve X$ the extension of $X$. We restrict ourselves to nonstandard arithmetic, i.e.\ extensions of the set $X=N$ of natural numbers (including zero). Objects of the form $\zve x$ for $x\in V(N)$ are called {\it standard}. For each $n\in N$, the element $\zve n$ is identified with $n$. Elements of $\zve N\setminus N$ are called {\it nonstandard integers}.\\ 

To every $k$-ary relation $\rho$ on $N$ corresponds a $k$-ary relation $\zve\rho$ on $\zve N$; the same holds for relations of a higher rank of the superstructure hierarchy (such as relations on subsets of $N$). To every operation $f:N^k\str N$ corresponds a $k$-ary operation $\zve f$ on $\zve N$. We will frequently use the extensions of addition, multiplication and the power operation, as well as the extensions of the usual orders $<$, $\leq$ and $\in$. To avoid overcomplicated formulas we will denote these operations and relations in the same way as their counterparts in $N$ (without a star). Moreover, we will assume that $\zvez\in$ is actually the membership relation on $V(\zve N)$ (for the justification of this see \cite{H}, Remark 5.1).

We learn about properties of such relations and functions mainly from the transfer principle. For example, we will use without mention facts such as $x^{y+1}=x^y\cdot x$ for $x,y\in\zve N$.\\

An object $x\in V(\zve N)$ is called {\it internal} if it is an element of a standard set $\zve A$ for $A\in V(N)\setminus N$. Thus all atoms $x\in\zve N$ are internal. Elements of internal sets are also internal.

\begin{pp}[The Internal Definition Principle]
For any formula $\varphi$, any $k\in\omega$ and any internal objects $a_1,a_2,\dots,a_n\in V(\zve N)$, the set 
$$\{x\in \zve V_k(N):\varphi(x,a_1,a_2,\dots,a_n)\mbox{ holds in }V(\zve N)\}$$
is internal.
\end{pp}

The important thing to remember is that the quantifiers in the Transfer Principle range only over internal objects. For example, by transfer it is easy to obtain the following.

\begin{pp}\label{intgreatest}
Every bounded internal subset of $\zve N$ has the greatest element.
\end{pp}

An important type of nonstandard extension is an {\it enlargement}. We call a binary relation $\rho$ in $V(N)$ {\it concurrent} if for every finitely many elements $a_1,a_2,\dots,a_k$ of its domain there is $b\in V(N)$ such that $a_i\rho b$ for $i=1,2,\dots,k$. A nonstandard extension $V(\zve N)$ is an enlargement if for every concurrent relation $\rho$ in $V(N)$ there is $x$ in the extension such that $\zve a\hspace{1mm}\zve\rho\hspace{1mm} x$ for all $a$ in the domain of $\rho$. Such extensions exist, see \cite{H}, Theorem 7.12.\\

{\bf The Stone-\v Cech compatification.} The set of all ultrafilters on $N$ is denoted by $\beta N$. For each $n\in N$ the principal ultrafilter $\{A\subseteq N:n\in A\}$ is identified with $n$. A topology can be defined on $\beta N$ so that it becomes the maximal compactification of the discrete space on $N$. This means that every $f:N\str N$ can be uniquely extended to a continuous function $\widetilde{f}:\beta N\str\beta N$. In this topology, for each $A\subseteq N$, $\overline{A}=\{\cF\in\beta N:A\in\cF\}$ is the closure of $A$. If an ultrafilter $\cF$ contains $A$ as an element, we will say that $\cF$ concentrates on $A$.

Only eventually constant sequences in $\beta N$ are convergent in the usual sense. Hence convergence via ultrafilters is often used: if $\cF,\cG_1,\cG_2,\dots$ are ultrafilters, $\lim_{n\str\cF}\cG_n=\cG$ if, for every $A\in\cG$, $\{n\in N:A\in\cG_n\}\in\cF$. More on these limits can be found in \cite{HS}, section 3.5.\\

For every $x\in\zve N$ the family $\{S\subseteq N:x\in\zve S\}$ is an ultrafilter; we denote this ultrafilter by $v(x)$. Thus a function $v:\zve N\str\beta N$ is obtained. For example, for $n\in N$ $v(n)$ is the corresponding principal ultrafilter.

In general, $v$ is not 1-1 (unless $V(\zve N)$ is obtained as ultrapower by a Hausdorff ultrafilter, see \cite{NF2}). $v$ is onto if $V(\zve N)$ is an enlargement.


\begin{fa}\label{rkpuritz}
(a) For every function $f:N\str N$ and every $x\in\zve N$, $\widetilde{f}(v(x))=v(\zve f(x))$.

(b) For every $f:N\str N$ and every $A\subseteq N$, $\zve(f[A])=\zve f[\zve A]$.
\end{fa}

\dokaz (a) is \cite{NR}, Lemma 1. (b) Let $B=f[A]$. Then, in $V(N)$, $(\zs n\in N)(n\in B\dl(\po m\in A)n=f(m))$. By transfer, for every $x\in\zve N$, $x\in\zve B\dl(\po y\in\zve A)x=\zve f(y)$.\kraj

Many aspects of the connection and the similarities between $\zve N$ and $\beta N$ were investigated in \cite{BNF} and \cite{NF1}.\\

The sets $\mu(\cF):=v^{-1}[\{\cF\}]$ for $\cF\in\beta N$ are called {\it monads}; they were investigated in a more general context in \cite{L}. By \cite{P2}, Theorem 3.1(a), monad of every nonprincipal ultrafilter in an enlargement has the same cardinality as $\zve N$ itself.\\

By \cite{P2}, Theorem 2.10, $\zve f[\mu(\cG)]\subseteq\mu(\widetilde{f}(\cG))$ for every $\cG\in\beta N$ and every $f:N\str N$. We make a short digression to provide more information about this in the following lemma.

\begin{lm}
(a) If $\mu(\cF)\neq\emptyset$, then the following conditions are equivalent: (i) $f[N]\in\cF$; (ii) $\mu(\cF)\subseteq\zve f[\zve N]$; (iii) $\mu(\cF)\cap\zve f[\zve N]\neq\emptyset$.

(b) If $\cF\in\beta N$ and $f:N\str N$ is such that $f[N]\in\cF$, then $\mu(\cF)=\bigcup\{\zve f[\mu(\cG)]:\widetilde{f}(\cG)=\cF\}$.
\end{lm}

\dokaz (a) If $f[N]\in\cF$, then for every $x\in\mu(\cF)$ we have $x\in\zve(f[N])=\zve f[\zve N]$ (Lemma \ref{rkpuritz}(b)), so $\mu(\cF)\subseteq\zve f[\zve N]$. (ii)$\Str$(iii) is obvious, and if $x\in\mu(\cF)\cap\zve f[\zve N]$, then $x\in\zve(f[N])$, so $f[N]\in\cF$.

(b) If $\widetilde{f}(\cG)=\cF$ then, by Lemma \ref{rkpuritz}(a), for every $x\in\mu(\cG)$ $v(\zve f(x))=\widetilde{f}(v(x))=\widetilde{f}(\cG)=\cF$, so $\zve f(x)\in\mu(\cF)$.

On the other hand, if $y\in\mu(\cF)$, by (a) the condition $f[N]\in\cF$ implies $y\in\zve f[\zve N]$, so there is $x\in\zve N$ such that $y=\zve f(x)$. If we denote $\cG=v(x)$, then $\widetilde{f}(\cG)=v(\zve f(x))=\cF$ and $y\in\zve f[\mu(\cG)]$.\kraj

{\bf Notation.} Throughout the paper, $N$ denotes the set of natural numbers (including zero) and $P$ denotes the set of (standard) prime numbers. 

If $n\in N$ is maximal such that $p^n\mid x$, we write $p^n\parallel x$. Analogous notation will be used in $\zve N$ (see Lemma \ref{transferi}(b)). For $z\in N$ let us denote $[0,z]_N=\{n\in N:n\leq z\}$. Analogously, $[0,z]_{\zve N}=\{n\in\zve N:n\leq z\}$ for $z\in\zve N$.

The elements of $\zve N$ will be denoted by small letters $x,y,z,\dots,$ with $p,q,\dots$ reserved for primes. The notation for ultrafilters will differ from that in \cite{So1} and \cite{So3}; they will be denoted by $\cF,\cG,\cH,\dots$, again with $\cP,\cQ,\dots$ reserved for prime ultrafilters.\\

For $A,B\subseteq N$ we denote $A\gstr=\{n\in N:\po a\in A\;a\mid n\}$, $A^2=\{a^2:a\in A\}$, $AB=\{ab:a\in A,b\in B,GCD(a,b)=1\}$ and $A^{(2)}=\{ab:a,b\in A,GCD(a,b)=1\}$. If $n\in N$, $nA=\{na:a\in A\}$. Also, $L_n=\{a_1a_2\dots a_n:a_1,a_2,\dots,a_n\in P\}$, $\cF\rest P=\{A\in\cF:A\subseteq P\}$ (for $\cF\in\beta N$) and $\cU=\{A\subseteq N:A\gstr=A\}$.\\

{\bf Divisibility in $\beta N$.} In \cite{So1} the author defined four relations on $\beta N$ extending divisibility in $N$. The one that most attention was given to is $\widemid$, further investigated in \cite{So3}:
$$\cF\widemid\cG\mbox{ iff }\cF\cap\cU\subseteq\cG.$$

We recapitulate some of the basic properties of this relation. It is not antisymmetric, so we think of it as an order on the equivalence classes $[\cF]$ of the relation defined by: $\cF=_\sim\cG\dl\cF\widemid\cG\land\cG\widemid\cF$. For such a class we denote $\mu([\cF])=\bigcup_{\cG=_\sim\cF}\mu(\cG)$.\\

An ultrafilter $\cP$ is {\it prime} (for $\widemid$) if it is divisible only by $1$ and itself. By \cite{So3}, Theorem 2.3, an ultrafilter $\cP$ is  prime if and only if $P\in\cP$. 

In \cite{So3} we described the lower part of the $\widemid$-hierarchy, more precisely the first $\omega$-many levels. $\widemid$ is antisymmetric within these lower levels (\cite{So3}, Lemma 5.13).

On the second level $\overline{L_2}$ (directly above prime ultrafilters) there are three types of ultrafilters: 

(1) those of the form $\cP^2$, generated by $\{A^2:A\in\cP\rest P\}$ for some prime $\cP\in\beta N$; 

(2) those containing $F_2^\cP=\{A^{(2)}:A\in\cP\rest P\}$ for some prime $\cP$ and 

(3) ultrafilters containing $F_{1,1}^{\cP,\cQ}=\{AB:A\in\cP\rest P,B\in\cQ\rest P,A\cap B=\emptyset\}$ for some two distinct prime ultrafilters $\cP$ and $\cQ$.

The ultrafilters of the third type are divisible by exactly two primes, and those of the first two types have only one prime divisor ("counted" twice). In a similar way, each ultrafilter on the $n$-th level $\overline{L_n}$ of the hierarchy has exactly $n$ "ingredients", not necessarilly distinct, with powers of primes $p^k$ counted $k$ times, see \cite{So3}, Theorem 5.5.

\section{The divisibility relation on $\zve N$}

In this section we recall some number-theoretic properties of the extension $\zvezmid$ of the divisibility relation $\mid$ on $N$. Note that the notions of greatest common divisor, least common multiplier and mutually prime numbers transfer directly from $N$ to $\zve N$. By transfer, $x\zvezmid y$ if and only if there is $k\in\zve N$ such that $y=kx$. If that is the case, we can write $\frac yx$ for $k$.

\begin{de}
$x\in\zve N\setminus\{0,1\}$ is prime if it is divisible only by $1$ and itself.
\end{de}

Clearly, "$x$ is prime" can be written as $(\zs y\in N)(y\mid x\Str y=1\lor y=x)$.

\begin{lm}\label{prosti}
For every $x\in\zve N$, $x$ is prime if and only if $x\in\zve P$.
\end{lm}

\dokaz The formula $(\zs x\in N\setminus\{0,1\})($"$x$ is prime"$\dl x\in P)$ holds in $V(N)$ so, by transfer, its counterpart $(\zs x\in\zve N\setminus\{0,1\})((\zs y\in\zve N)(y\mid x\Str y=1\lor y=x)\dl x\in\zve P)$ holds in $V(\zve N)$.\kraj

The next lemma also follows directly from the Transfer Principle.

\begin{lm}\label{transferi3}
(a) For all $z\in\zve N$, $n\in N$ and $A\subseteq N$: $z\in\zve(nA)$ if and only if $z=nx$ for some $x\in\zve A$.

(b) For all $x\in\zve N$, $n\in N$ and $A\subseteq N$: $nx\in\zve(nA)$ if and only if $x\in\zve A$.
\end{lm}

In particular, $z\in\zve N$ is divisible by $n\in N$ if and only if $z\in\zve(nN)$. Our next lemma lists several other properties of $\zvezmid$ that mostly follow directly from the Transfer Principle.

\begin{lm}\label{transferi}
(a) $|\zve P|=|\zve N|$.

(b) For every $x\in\zve N\setminus\{0\}$ and every $p\in\zve P$ there is maximal $a\in\zve N$ such that $p^a\zvezmid x$.

(c) If $x,y\in\zve N$ have the same sets of divisors of the form $p^z$ ($p\in\zve P$, $z\in\zve N$), then $x=y$.
\end{lm}

\dokaz (a) Let $p:N\str P$ be the function mapping every $n\in N$ to the $n$-th prime number; so $p(0)=2$, $p(1)=3$ etc. Then the formula $(\zs m\in N)$"$m$ is prime"$\dl(\po n\in N)m=p(n)$ holds in $V(N)$, so by transfer each $\zve p(x)$ is prime in $V(\zve N)$. $p$ is a bijection, so $\zve p$ is a bijection too. Finally, by Lemma \ref{rkpuritz}(b), $\zve P=\zve p[\zve N]$.

(b) We have $(\zs x\in N\setminus\{0\})(\zs p\in P)(\po a\in N)(p^a\mid x\land p^{a+1}\nmid x)$, so the same holds in $V(\zve N)$.

(c) Since $(\zs x,y\in N)((\zs p\in P)(\zs n\in N)(p^n\mid x\dl p^n\mid y)\Str x=y)$ holds, it follows that for all $x,y\in\zve N$: $(\zs p\in\zve P)(\zs z\in\zve N)(p^z\zvezmid x\dl p^z\zvezmid y)$ implies $x=y$.\kraj

$\zve N$ is not well-ordered, so infinite sums and products can not be defined in the usual way, by induction. However, using transfer we can bypass this, using an idea described in \cite{H}, Remark 5.8.

The next theorem is an extension of the fundamental theorem of arithmetic. Within it $p$ is the function enumerating all primes (defined in the proof of Lemma \ref{transferi}(a)). Note that, since a sequence $f:[0,z]_N\str N$ is a set of ordered pairs, we have $f\in V_3(N)$.

\begin{te}\label{fund}
(a) For every $z\in\zve N$ and every internal sequence $\langle h(n):n\leq z\rangle$ there is unique $x\in\zve N$ such that $\zve p(n)^{h(n)}\zvepar x$ for $n\leq z$ and $\zve p(n)\nzvemid x$ for $n>z$; we denote such element by $\prod_{n\leq z}\zve p(n)^{h(n)}$.

(b) Every $x\in\zve N$ can be uniquely represented as $\prod_{n\leq z}\zve p(n)^{h(n)}$ for some $z\in\zve N$ and some internal sequence $\langle h(n):n\leq z\rangle$ such that $h(z)>0$.
\end{te}

\dokaz (a) Let, as usual, "$f:X\str N$" denote the formula: $f\mbox{ is a function}\land\dom(f)=X\land\ran(f)\subseteq N$. Let "$g(z)=\prod_{n\leq z}f(n)$" denote the formula 
$$g(0)=f(0)\land(\zs n<z)g(n+1)=g(n)\cdot f(n+1).$$
In $V(N)$ we have:

\parbox{10cm}{\begin{eqnarray*}
& (\zs z\in N)(\zs f\in V_3(N))(f:[0,z]_N\str N\Str(\po_1 g\in V_3(N))\nonumber\\
& (g:[0,z]_N\str N\land g(z)=\prod_{n\leq z}f(n)).
\end{eqnarray*}}\hfill\parbox{1cm}{\begin{eqnarray}\label{proizvod}\end{eqnarray}}\\

By transfer, the same holds in $V(\zve N)$. Now, if we are given $z\in\zve N$ and an internal sequence $\langle h(n):n\leq z\rangle$, by the Internal Definition Principle the sequence $\langle f(n):n\leq z\rangle$ defined by $f(n)=\zve p(n)^{h(n)}$ is also internal. (\ref{proizvod}) now produces a sequence $\langle g(n):n\leq z\rangle$ such that $g:[0,z]_{\zve N}\str\zve N\land g(z)=\prod_{n\leq z}\zve p(n)^{h(n)}$.

Now we use transfer again: in $V(N)$
\begin{eqnarray*}
&(\zs z\in N)(\zs h\in V_3(N))(\zs g\in V_3(N))(h:[0,z]_N\str N\land\\
&g:[0,z]_N\str N\land g(z)=\prod_{n\leq z}p(n)^{h(n)}\Str\\
&(\zs n\leq z)p(n)^{h(n)}\parallel g(z)\land(\zs n>z)p(n)\nmid g(z)),
\end{eqnarray*}
so the same holds in $V(\zve N)$ and $x=g(z)$ is the wanted element. Uniqueness follows from Lemma \ref{transferi}(c).

(b) For $x\in\zve N$, the set of primes that divide $x$ is clearly bounded ($x$ can not be divisible by primes greater than itself). It is also internal by the Internal Definition Principle, so it has the greatest element by Proposition \ref{intgreatest}; let $z$ be this element. For $n\leq z$ we define $h(n)$ to be the greatest $a\in\zve N$ such that $\zve p(n)^a\;\zvezmid x$; such $a$ exists by Lemma \ref{transferi}(b) and the obtained sequence $\{(n,a):n\leq z\land\zve p(n)^a\zvepar x\}$ is internal, again by the Internal Definition Principle. Now by (a) we get $x':=\prod_{n\leq z}\zve p(n)^{h(n)}$ divisible by the same powers of primes as $x$. By Lemma \ref{transferi}(c), $x=x'$.

To prove uniqueness, assume $x=\prod_{n\leq z'}\zve p(n)^{h'(n)}$ for some $z'\in\zve N$ and some sequence $\langle h'(n):n\leq z'\rangle$. If $z'>z$, this would mean that $x$ is divisible by $p(z')$; if $z'<z$ then $x$ would not be divisible by $p(z)$. Either way we reach a contradiction, so $z'=z$. In a similar manner we get a contradiction if we assume that $h'(n)\neq h(n)$ for some $n\leq z$.\kraj

Let $\lev:N\setminus\{0\}\str N$ be the function calculating the level of each $n\in N\setminus\{0\}$ in the $\mid$-hierarchy. More precisely, if $n=p_1^{a_1}p_2^{a_2}\dots p_k^{a_k}$, let $\lev(n)=a_1+a_2+\dots+a_k$. Its extension $\zve\lev$ does the same for elements in $\zve N$, when represented as $\prod_{n\leq z}\zve p(n)^{h(n)}$, as in Theorem \ref{fund}. Namely, let "$g(z)=\sum_{n\leq z}f(n)$" denote the formula 
$$g(0)=f(0)\land(\zs n<z)g(n+1)=g(n)+f(n+1).$$
Then, by transfer, $\zve\lev$ is the unique function satisfying, for every $z\in\zve N$ and every internal $h:[0,z]\str N$, the formula $g(z)=\prod_{n\leq z}\zve p(n)^{h(n)}\land\zve\lev(g(z))=\sum_{n\leq z}h(n)$.

Here are some properties of the function $\zve\lev$, proven easily by transfer.

\begin{lm}\label{levels}
Let $x,y\in\zve N$ be such that $x\zvezmid y$.

(a) Either $x=y$ or $\zve\lev(x)<\zve\lev(y)$.

(b) If $\zve\lev(x)<a<\zve\lev(y)$, then there is $z\in\zve N$ such that $\zve\lev(z)=a$, $x\zvezmid z$ and $z\zvezmid y$.
\end{lm}

In the following two sections we will see that the nice structure of the $\zvezmid$-hierarchy is mostly transferred to the first $\omega$-many levels of the $\widemid$, but not above them.

\section{The connection with the Stone-\v Cech compactification}

We have already encountered several analogies between the divisibility relation $\zvezmid$ on $\zve N$ and the relation $\widemid$ on $\beta N$. First, $x\in\zve N$ is divisible by $n\in N$ if and only if $x\in\zve(nN)$ (Lemma \ref{transferi3}(a)) and $\cF\in\beta N$ is divisible by $n$ if and only if $nN\in\cF$ (\cite{So1}, Lemma 5.1). Also, an ultrafilter $\cP$ is prime if and only if it concentrates on the set of primes. By Lemma \ref{prosti} the same thing holds in $\zve N$ for the relation $\zvezmid$, so $x\in\zve N$ is prime if and only if $v(x)$ is a prime ultrafilter.

We will now establish a connection between the relation $\widemid$ and the divisibility in $\zve N$, showing that these similarities are not coincidental. It also shows that $\widemid$ is, in some sense, "the right" divisibility relation to investigate in $\beta N$.

\begin{te}\label{ekviv}
The following conditions are equivalent for every two ultrafilters $\cF,\cG\in\beta N$:

(i) $\cF\widemid\cG$;

(ii) in every enlargement $V(\zve N)$, there are $x,y\in\zve N$ such that $v(x)=\cF$, $v(y)=\cG$ and $x\zvezmid y$;

(iii) in some enlargement $V(\zve N)$, there are $x,y\in\zve N$ such that $v(x)=\cF$, $v(y)=\cG$ and $x\zvezmid y$.
\end{te}

\dokaz (i)$\Str$(ii) Let $V(\zve N)$ be an enlargement, and let $\cF,\cG\in\beta N$ be such that $\cF\widemid\cG$. We define a binary relation $\rho\subseteq P(N)\times N^2$:
$$A\rho(m,n)\mbox{ iff }(m\in A\dl A\in\cF)\land(n\in A\dl A\in\cG)\land m\mid n.$$
We prove that $\rho$ is concurrent. Let a finite number of subsets of $N$ be given; we need to find a pair $(m,n)\in N^2$ such that $A\rho(m,n)$ for all given sets $A$. Since ultrafilters (and their complements) are closed for finite intersections, we may assume that we have at most four given sets: $A_1\in\cF\cap\cG$, $A_2\in\cF\setminus\cG$, $A_3\in\cG\setminus\cF$ and $A_4\notin\cF\cup\cG$ (if there are more sets from the same class, say $A_1^1,A_2^1,\dots,A_k^1\in\cF\cap\cG$, we can replace them with their intersection). So we are looking for a pair $(m,n)\in N^2$ such that $m\in A_1\cap A_2\setminus(A_3\cup A_4)$, $n\in A_1\cap A_3\setminus(A_2\cup A_4)$ and $m\mid n$. The set $B:=A_1\cap A_2\setminus(A_3\cup A_4)$ belongs to $\cF$. $B\gstr=\{n\in N:\po b\in B\;b\mid n\}$ belongs to $\cF\cap\cU$, so it must be in $\cG$ as well (since $\cF\widemid\cG$). Since the set $C:=A_1\cap A_3\setminus(A_2\cup A_4)$ is also in $\cG$, we can choose $n\in C\cap B\gstr$. But $n\in B\gstr$ means that there is $m\in B$ such that $m\mid n$.

Now, since $\zve N$ is an enlargement, there is a pair $(x,y)\in\zve(N^2)=(\zve N)^2$ such that $\zve A\hspace{1mm}\zve\rho\hspace{1mm}(x,y)$ for all $A\in P(N)$. By transfer, $A\in\cF$ holds if and only if $\zve A\in\zve\cF$. Thus we get, for all $A\in P(N)$: 
$$(x\in\zve A\dl A\in\cF)\land(y\in\zve A\dl A\in\cG)\land x\zvezmid y.$$
This means that $v(x)=\cF$ and $v(y)=\cG$.

(ii)$\Str$(iii) is obvious.

(iii)$\Str$(i) In $V(N)$ we have, for every $A\in\cU$: $(\zs m\in A)(\zs n\in N)(m\mid n\Str n\in A)$. Hence the same holds for $\zve A$ in any extension $V(\zve N)$ i.e.\ $\zve A$ is closed upwards for $\zvezmid$. This means that, if $x\zvezmid y$ for some $x\in\mu(\cF)$, $y\in\mu(\cG)$, then $x\in\zve A$ implies $y\in\zve A$ for every $A\in\cU$. Thus $\cF\cap\cU\ps\cG$, i.e.\ $\cF\widemid\cG$.\kraj

Note that the implication (iii)$\Str$(i) holds in any extension, not only in an enlargement.

The next example shows that whether or not $x\zvezmid y$ holds is not independent from the choice of representatives $x\in\mu(\cF)$ and $y\in\mu(\cG)$.

\begin{ex}
Since $\zvezmid$ is reflexive (by transfer), it suffices to find $a,b\in\zve N$ such that $v(a)=v(b)$ and $a\nzvemid b$. So assume $x,y\in\zve N$ are such that $v(x)=v(y)$, $x\neq y$ and $x\zvezmid y$. Let $f:N\str N$ be a function such that
\begin{equation}\label{nedeli}
(\forall m,n\in N)(m\mid n\land m\neq n\Str f(m)\nmid f(n))
\end{equation}
($f(n)$ is easily constructed by recursion on $n\in N$). Now let $a=\zve f(x)$ and $b=\zve f(y)$. Then $a\zvezmid a$, 
$$v(a)=v(\zve f(x))=\widetilde{f}(v(x))=\widetilde{f}(v(y))=v(\zve f(y))=v(b)$$
(by Fact \ref{rkpuritz}(a)), but $a\nzvemid b$, by (\ref{nedeli}) and transfer.
\end{ex}

A chain (in $(N,\mid)$) is a set $C\ps N$ of elements linearly ordered by $\mid$; an antichain is a set $A\ps N$ of $\mid$-incomparable elements. The family of subsets of $N$ containing the complements of all chains and the complements of all antichains has the finite intersection property, so there are ultrafilters containing no chains and no antichains. On the other hand, every selective ultrafilter contains at least one of these two types of sets.

\begin{lm}\label{achains}
Let $V(\zve N)$ be an enlargement and $\cF\in\beta N\setminus N$. 

(a) $\cF$ contains no infinite antichains as elements if and only if there are distinct $x,y\in\mu(\cF)$ such that $x\zvezmid y$.

(b) $\cF$ contains no infinite chains as elements if and only if there are distinct $x,y\in\mu(\cF)$ such that neither $x\zvezmid y$ nor $y\zvezmid x$.
\end{lm}

\dokaz (a) First assume there is an infinite antichain $A\in\cF$. In $V(N)$ we have $(\zs m,n\in A)m\nmid n$, so the same holds in $V(\zve N)$. Thus there are no $\zvezmid$-comparable elements in $\zve A$, so there are none in $\mu(\cF)$.

Now let $\cF$ contain no infinite antichains. We define a binary relation $\rho\ps P(N)\times N^2$:
$$A\rho(m,n)\mbox{ iff }(m\in A\dl A\in\cF)\land(n\in A\dl A\in\cF)\land m\neq n\land m\mid n.$$
$\rho$ is concurrent: if we are given finitely many subsets of $N$, let $A_1$ be the intersection of those in $\cF$, and $A_2$ the intersection of those outside $\cF$. Then $A_1\setminus A_2\in\cF$, so it is not an antichain. Hence there are distinct $m,n\in A_1\setminus A_2$ such that $m\mid n$.

$V(\zve N)$ is an enlargement, so there are distinct $x,y\in\zve N$ such that $x\zvezmid y$ and $x,y\in\zve A$ for all $A\in\cF$; but then $x,y\in\mu(\cF)$.

The proof for (b) is analogous.\kraj

Since ultrafilters on the first $\omega$-many levels of $\widemid$-hierarchy contain antichains (an ultrafilter on the $n$-th level contains the set $L_n$), using Lemmas \ref{levels} and \ref{achains}(a) by induction on $n$ we can easily prove that, if $\cF\in\overline{L_n}$ then $\zve\lev(x)=n$ for every $x\in\mu(\cF)$. Hence such ultrafilters correspond precisely to the first $\omega$-many levels of $\zvezmid$-hierarchy in an enlargement (containing $x\in\zve N$ such that $\zve\lev(x)$ is finite). We will investigate in more detail the connection of $x\in\zve N$ such that $\zve\lev(x)=2$ with their corresponding ultrafilters $v(x)$ described in the Introduction, obtaining in particular a better insight into the origin of the ultrafilters containing $F_2^\cP$ (type (2) from the Introduction). 

\begin{lm}
Let $V(\zve N)$ be any nonstandard extension.

(a) $x\in\zve N$ is of the form $p^2$ for some $p\in\zve P$ if and only if $v(x)=\cP^2$ for some prime ultrafilter $\cP$.

(b) $x\in\zve N$ is of the form $p\cdot q$ for two distinct primes $p,q$ such that $v(p)=v(q)=\cP$ if and only if $v(x)\supseteq F_2^\cP$.

(c) $x\in\zve N$ is of the form $p\cdot q$ for two primes $p,q$ such that $v(p)=\cP$, $v(q)=\cQ$ and $\cP\neq\cQ$ if and only if $v(x)\supseteq F_{1,1}^{\cP,\cQ}$.
\end{lm}

\dokaz (a) Let $sq:N\str N^2$ be the squaring function: $sq(n)=n^2$ for $n\in N$. Then $x=\zve sq(p)$ for some $p\in\zve P$ implies $v(x)=\widetilde{sq}(v(p))=(v(p))^2$ and $v(p)\in\overline{P}$.

Now let $v(x)=\cP^2$ for a prime ultrafilter $\cP$. Then $x\in\zve(P^2)=\zve sq[\zve P]$ (by Lemma \ref{rkpuritz}(b)) so $x=\zve sq(p)=p^2$ for some $p\in\zve P$.

(b) Let $x=p\cdot q$ for some $p,q\in\zve P$ such that $v(p)=v(q)=\cP$. Let $A\in\cP\rest P$. Since $(\zs a,b\in A)(a\neq b\Str ab\in A^{(2)})$, by transfer 
\begin{equation}\label{transfer2}
(\zs a,b\in\zve A)(a\neq b\Str ab\in\zve(A^{(2)}))
\end{equation}
so $x=pq\in\zve(A^{(2)})$. Thus $A^{(2)}\in v(x)$.

For the other direction let $x\in\zve N$ be such that $v(x)\supseteq F_2^\cP$ for some prime ultrafilter $\cP$. Then $x\in\zve(P^{(2)})$. In $V(N)$ we have $(\zs n\in P^{(2)})(\po a,b\in P)(a\neq b\land n=ab)$ so, by transfer, $x=pq$ for some distinct $p,q\in\zve P$. To prove that $p,q\in\zve A$ for each $A\in\cP\rest P$, assume the opposite: either $p,q\in\zve(P\setminus A)$ or one of them belongs to $\zve A$ and the other to $\zve(P\setminus A)$. But if, for example, the first option holds then, as in (\ref{transfer2}), we get $pq\in\zve((P\setminus A)^{(2)})\subseteq\zve(P^{(2)}\setminus A^{(2)})$, a contradiction since $A^{(2)}\in F_2^\cP$.

(c) The proof is similar to the proof of (b).\kraj

Thus, ultrafilters containing families of the form $F_2^\cP$ actually have two distinct "ingredients", but such that $\beta N$ can not distinguish between.


\section{Above finite levels}

There are, of course, also ultrafilters not concetrating on any $L_n$ for $n\in N$. The investigation of $\widemid$ becomes much more complicated at these higher levels. Limits of ultrafilters will prove useful for this purpose.

\begin{lm}\label{supremum}
Let $\langle \cG_n:n\in N\rangle$ be a $\widemid$-increasing sequence in $\beta N$.

(a) $\lim_{n\str\cF}\cG_n\cap\cU=\bigcup_{n\in N}(\cG_n\cap\cU)$ for any nonprincipal ultrafilter $\cF$.

(b) For any two nonprincipal ultrafilters $\cF_1$ and $\cF_2$, 
$$\lim_{n\str\cF_1}\cG_n=_\sim\lim_{n\str\cF_2}\cG_n.$$

(c) Let $\cG=\lim_{n\str\cF}\cG_n$ for some $\cF$. If $\cW\in\beta N$ is such that $\cG_n\widemid\cW$ for all $n\in N$, then $\cG\widemid\cW$.
\end{lm}

\dokaz (a) If $A\in\cG_m\cap\cU$ for some $m\in N$, then $A\in\cG_n\cap\cU$ for all $n\geq m$. Hence the set $\{n\in N:A\in\cG_n\}$ is cofinite, so it belongs to $\cF$. It follows that $A\in\lim_{n\str\cF}\cG_n$.

On the other hand, assume $A\in\cU$ is such that $A\notin\cG_n$ for all $n\in N$. Then $N\setminus A\in\cG_n$ for all $n$, so $N\setminus A\in\lim_{n\str\cF}\cG_n$ and $A\notin\lim_{n\str\cF}\cG_n$.

(b) Follows from (a).

(c) $\cG_n\widemid\cW$ means that $\cG_n\cap\cU\subseteq\cW$. If this holds for all $n\in N$, by (a) $\cG\cap\cU\subseteq\cW$, so $\cG\widemid\cW$.\kraj

In view of Lemma \ref{supremum}(b), we will write $[\cG]=\lim_{n\str\infty}\cG_n$ if $\cG=\lim_{n\str\cF}\cG_n$ for some nonprincipal $\cF$. 

\begin{ex}
There are also infinite $\widemid$-decreasing sequences in $\beta N$. Let $\{p_i:i\in N\}$ be an enumeration of $P$. Let $r_m^n=\prod_{i=m}^np_i$, and $[\cG_m]=\lim_{n\str\infty}r_m^{m+n}$. Then each $\cG_m$ is divisible by all $r_n^{n+k}$ for $n\geq m$ so, by Lemma \ref{supremum}(c), $\cG_n\widemid\cG_m$ for $n>m$. Also, if $n>m$ then $p_m\widemid\cG_m$ but $p_m\nwidemid\cG_n$, so the sequence $\langle \cG_n:n\in N\rangle$ is strictly decreasing.
\end{ex}

The following lemma is proved analogously to Lemma \ref{supremum}.

\begin{lm}
Let $\langle \cG_n:n\in N\rangle$ be a $\widemid$-decreasing sequence in $\beta N$. 

(a) $\lim_{n\str\cF}\cG_n\cap\cU=\bigcap_{n\in N}(\cG_n\cap\cU)$ for any nonprincipal ultrafilter $\cF$.

(b) For any two nonprincipal ultrafilters $\cF_1$ and $\cF_2$,
$$\lim_{n\str\cF_1}\cG_n=_\sim\lim_{n\str\cF_2}\cG_n.$$

(c) Let $\cG=\lim_{n\str\cF}\cG_n$ for some $\cF$. If $\cW\in\beta N$ is such that $\cW\widemid\cG_n$ for all $n\in N$, then $\cW\widemid\cG$.
\end{lm}

\begin{lm}\label{sled}
Let $x\in\zve N$, $p\in P$, $\cF=v(x)$ and $\cG=v(px)$. If there is $n\in N$ such that $p^n\zvepar x$, then $\cG$ is an immediate successor of $\cF$ in $(\beta N,\widemid)$.
\end{lm}

\dokaz By Theorem \ref{ekviv}, $\cF\widemid\cG$. Since $p^{n+1}N\in(\cG\cap\cU)\setminus\cF$, $\cF\neq_\sim\cG$. It remains to show that there is no $\cH\in\beta N$ such that $\cF\widemid\cH$, $\cH\widemid\cG$ but $\cF\neq_\sim\cH\neq_\sim\cG$.

{\it Claim.} The set $\cG\cap\cU$ is generated by $(\cF\cap\cU)\cup\{p^{n+1}N\}$. This means that every ultrafilter containing $(\cF\cap\cU)\cup\{p^{n+1}N\}$ is divisible by $\cG$. So assume $B\in\cG\cap\cU$ and let $B':=B\cap(p^{n+1}N\setminus p^{n+2}N)$. Clearly $B'\in\cG$. We can write $B'$ in the form $B'=p^{n+1}A=p(p^nA)$ for some $A$ such that $A\cap pN=\emptyset$. By Lemma \ref{transferi3}(b), $p^nA\in\cF$. But $(p^{n+1}A)\gstr=(p^nA)\gstr\cap p^{n+1}N$ and $B\supseteq(B')\gstr=(p^{n+1}A)\gstr$, so $B$ must belong to every ultrafilter containing $(\cF\cap\cU)\cup\{p^{n+1}N\}$.

Now assume $\cH$ is as above. If $p^{n+1}N\in\cH$, by Claim we have $\cH=_\sim\cG$. Otherwise, if we assume $A\in(\cH\cap\cU)\setminus\cF$, then $A_1:=A\cap(p^nN\setminus p^{n+1}N)\in\cH\setminus\cF$. But $A_2:=(p^nN\setminus p^{n+1}N)\setminus A_1\in\cF$ is disjoint from $A_1$, so $A_2\gstr\in(\cF\cap\cU)\setminus\cH$. This contradicts $\cF\widemid\cH$.\kraj



The following example shows that the condition of existence of $n\in N$ such that $p^n\zvepar x$ can not be eliminated from the lemma above.

\begin{ex}\label{exlimit}
Let $[\cG]=\lim_{n\str\infty}p^n$ for some $p\in P$. We show that $\mu([\cG])=\{p^x:x\in\zve N\setminus N\}$. By Lemma \ref{supremum}(a) $\cG\cap\cU=\{A\in\cU:p^nN\subseteq A\mbox{ for some }n\in N\}=\{A\in\cU:p^n\in A\mbox{ for some }n\in N\}$.

First, no elements divisible by any prime other than $p$ can belong to $\mu([\cG])$: if $x\in\mu([\cG])$ and $q\in\zve P\setminus\{p\}$ are such that $q\zvezmid x$, there is $A\subseteq P$ such that $p\notin A$ but $q\in\zve A$, so $A\gstr\in(v(x)\cap\cU)\setminus\cG$, a contradiction.

Now we prove that $v(p^x)=_\sim\cG$ for $x\in\zve N\setminus N$. Assume the opposite, that there is $A\in\cU\setminus\cG$ such that $p^x\in\zve A$ for some $x\in\zve N$. Then $A\cap\{p^n:n\in N\}=\emptyset$ (otherwise $A\in\cG$). By the Transfer Principle $\zve A\cap\{p^x:x\in\zve N\}=\emptyset$, a contradiction.
\end{ex}

It is easy to see that there is the $\widemid$-maximal class in $\beta N$ (since $\cU$ has the finite intersection property). Let $MAX$ denote this maximal class.

\begin{lm}\label{maks}
For every $x\in\zve N$: $x\in\mu(MAX)$ if and only if $n\zvezmid x$ for all $n\in N$.
\end{lm}

\dokaz If $x\in\mu(MAX)$, then $x\in\zve A$ for all $A\in\cU$ so, for every $n\in N$, $x\in\zve(nN)$, which is equivalent to $n\zvezmid x$ by Lemma \ref{transferi3}(a).

Now assume $n\zvezmid x$ (i.e.\ $nN\in v(x)$) for all $n\in N$. Let $A\in\cU$ be arbitrary and let $m$ be a $\mid$-minimal element of $A$. Then $mN\subseteq A$, so $A\in v(x)$. Hence $\cU\subseteq v(x)$, which means that $x\in\mu(MAX)$.\kraj

In particular, this means that $MAX=\lim_{n\str\infty}n!$. Thus the distribution of ultrafilters by levels, described in previous section, fails for ultrafilters above finite levels by Lemma \ref{supremum}(c): $MAX$ would be on the $\omega$-th level, and at the same time has predecessors on infinite levels.

\begin{te}\label{successor}
(a) Every $\cF\in\beta N\setminus MAX$ has an immediate successor in $(\beta N,\widemid)$.

(b) Every $\cF\in\beta N$ such that there are $p\in P$ and $n\in N\setminus\{0\}$ so that $p^n\zvepar\cF$ has an immediate predecessor in $(\beta N,\widemid)$.
\end{te}

\dokaz (a) Let $x\in\mu(\cF)$. Since $x\notin\mu(MAX)$, by Lemma \ref{maks} there is $m\in N$ such that $m\nmid x$. Thus there are $p\in P$ and $n\in N$ so that $p^n\zvepar x$. By Lemma \ref{sled} $v(px)$ is the immediate successor of $\cF$.

(b) Let $x\in\mu(\cF)$. As in (a) we can show that $v\left(\frac xp\right)$ is the immediate predecessor of $\cF$.\kraj

\begin{ex}
We show that the condition $n\neq 0$ of Theorem \ref{successor}(b) can not be omitted. Let $[\cG]=\lim_{n\str\infty}p^n$, as in Example \ref{exlimit}. Then $\cG$ is divisible by all powers of $p$ and not divisible by any other prime. Assume $\cG$ has an immediate predecessor $\cF$. We consider two cases. 

$1^\circ$ $p^n\widemid\cF$ for all $n\in N$. Then by Lemma \ref{supremum}(c) $\cG\widemid\cF$, a contradiction. 

$2^\circ$ $p^n\hspace{1mm}\widetilde{\parallel}\hspace{1mm}\cF$ for some $n$. Let $x\in\mu(\cF)$ and $y\in\mu(\cG)$ be such that $x\zvezmid y$. Then, by Lemma \ref{sled}, $v(px)$ is a successor of $\cF$ and, since $p^{n+1}\zvezmid y$ and $px$ is the least common multiplier of $x$ and $p^{n+1}$, $px\zvezmid y$ as well, meaning that $v(px)\widemid\cG$. In a similar way we obtain $v(p^2x)\widemid\cG$ so, since $v(px)\neq v(p^2x)$, $\cF$ is not an immediate predecessor of $\cG$.
\end{ex}

\section{Open problems and final remarks}

We mention several questions, the answers to which may shed some more light to the above results.

\begin{qu}\label{que1}
(a) Does $v(x)=v(y)$ imply $v(xz)=v(yz)$ for $x,y,z\in\zve N$? 

(b) More generally, does $v(x)=v(y)$ and $v(z)=v(u)$ imply $v(xz)=v(yu)$ for $x,y,z,u\in\zve N$?
\end{qu}

\begin{qu}
Does $v(x)=_\sim v(y)$ imply $v(xz)=_\sim v(yz)$ for $x,y,z\in\zve N$?
\end{qu}

This is true if $x,y,z\in\zve P$ are distinct. Namely, assume there is $A\in(v(xz)\cap\cU)\setminus v(yz)$. Let $X\in v(x)$, $Y\in v(y)$ and $Z\in v(z)$ be subsets of $P$ such that $X,Y$ and $Z$ are disjoint. By \cite{So3} Lemma 3.7, $XZ\in v(xz)$ and $YZ\in v(yz)$. Then $A':=A\cap XZ\in v(xz)$. If we define $f:P^{(2)}\str P$ by $f(ab)=a$ (for $a\in X\cup Y,b\in Z$) and $f(n)$ arbitrary if $n\notin (X\cup Y)Z$, then $\widetilde{f}(v(xz))=v(x)$, $\widetilde{f}(v(yz))=v(y)$, $f[A']\subseteq X$ and $f[A']\gstr\in(v(x)\cap\cU)\setminus v(y)$, a contradiction with $v(x)=_\sim v(y)$.

\begin{qu}
Let $\zve N$ be an enlargement.

(a) Does $\cF\widemid\cG$ imply $(\zs x\in\mu(\cF))(\po y\in\mu(\cG))x\zvezmid y$?

(b) Does $\cF\widemid\cG$ imply $(\zs y\in\mu(\cG))(\po x\in\mu(\cF))x\zvezmid y$?
\end{qu}

Part (a) is true if the answer to Question \ref{que1} is "yes". Namely, $\cF\widemid\cG$ means that there are $u\in\mu(\cF)$ and $z\in\zve N$ such that $uz\in\mu(\cG)$. But then it would follow that, for every $x\in\mu(\cF)$, $xz\in\mu(\cG)$ as well.

\begin{qu}
Let us call a set $X\subseteq\zve N$ convex if for all $x,y\in X$ and $z\in\zve N$, $x\zvezmid z$ and $z\zvezmid y$ implies $z\in X$. Is $\mu(\cF)$ a convex set for every $\cF\in\beta N$?
\end{qu}

Clearly, every $\mu([\cF])$ is convex: if $x\in\mu(\cG_1)$ and $y\in\mu(\cG_2)$, then $\cG_1=_\sim\cG_2=_\sim\cF$ would imply that $\cG_1\cap\cU=\cG_2\cap\cU$, so $x\zvezmid z$ and $z\zvezmid y$ would imply $\cG_1\cap\cU\subseteq v(z)\cap\cU\subseteq\cG_2\cap\cU$, so $v(z)=_\sim\cF$.

\begin{qu}
Can we strengthen Example \ref{exlimit} in the following sense: if $[\cF]=\lim_{n\str\infty}v(xp^n)$ for some $x\in\zve N$, does $\mu([\cF])=\{xp^n:n\in\zve N\setminus N\}$?
\end{qu}

The research was supported by the Ministry of Education, Science and Technological Development of the Republic of Serbia (project 174006).

\end{document}